# Comparison of robust tests for genetic association using case-control studies

**Gang Zheng[1], Boris Freidlin[2] and Joseph L. Gastwirth[3,*]**

*National Heart, Lung and Blood Institute, National Cancer Institute and George Washington University*

**Abstract:** In genetic studies of complex diseases, the underlying mode of inheritance is often not known. Thus, the most powerful test or other optimal procedure for one model, e.g. recessive, may be quite inefficient if another model, e.g. dominant, describes the inheritance process. Rather than choose among the procedures that are optimal for a particular model, it is preferable to see a method that has high efficiency across a family of scientifically realistic models. Statisticians well recognize that this situation is analogous to the selection of an estimator of location when the form of the underlying distribution is not known. We review how the concepts and techniques in the efficiency robustness literature that are used to obtain efficiency robust estimators and rank tests can be adapted for the analysis of genetic data. In particular, several statistics have been used to test for a genetic association between a disease and a candidate allele or marker allele from data collected in case-control studies. Each of them is optimal for a specific inheritance model and we describe and compare several robust methods. The most suitable robust test depends somewhat on the range of plausible genetic models. When little is known about the inheritance process, the maximum of the optimal statistics for the extreme models and an intermediate one is usually the preferred choice. Sometimes one can eliminate a mode of inheritance, e.g. from prior studies of family pedigrees one may know whether the disease skips generations or not. If it does, the disease is much more likely to follow a recessive model than a dominant one. In that case, a simpler linear combination of the optimal tests for the extreme models can be a robust choice.

## 1. Introduction

For hypothesis testing problems when the model generating the data is known, optimal test statistics can be derived. In practice, however, the precise form of the underlying model is often unknown. Based on prior scientific knowledge a family of possible models is often available. For each model in the family an optimal test statistic is obtained. Hence, we have a collection of optimal test statistics corresponding to each member of the family of scientifically plausible models and need to select one statistic from them or create a robust one, that combines them. Since using any single optimal test in the collection typically results in a substantial loss of efficiency or power when another model is the true one, a robust procedure with reasonable power over the entire family is preferable in practice.

*Supported in part by NSF grant SES-0317956.

[1]Office of Biostatistics Research, National Heart, Lung and Blood Institute, Bethesda, MD 20892-7938, e-mail: zhengg@nhlbi.nih.gov

[2]Biometric Research Branch, National Cancer Institute, Bethesda, MD 20892-7434, e-mail: freidlinb@ctep.nci.nih.gov

[3]Department of Statistics, George Washington University, Washington, DC 20052, e-mail: jlgast@gwu.edu

The above situation occurs in many applications. For example, in survival analysis Harrington and Fleming [14] introduced a family of statistics $G^\rho$. The family includes the log-rank test ($\rho = 0$) that is optimal under the proportional hazards model and the Peto-Peto test ($\rho = 1$, corresponding to the Wilcoxon test without censoring) that is optimal under a logistic shift model. In practice, when the model is unknown, one may apply both tests to survival data. It is difficult to draw scientific conclusions when one of the tests is significant and the other is not. Choosing the significant test *after* one has applied both tests to the data increases the Type I error. A second example is testing for an association between a disease and a risk factor in contingency tables. If the risk factor is a categorical variable and has a natural order, e.g., number of packs of cigarette smoking per day, the Cochran-Armitage trend test is typically used. (Cochran [3] and Armitage [1]) To apply such a trend test, increasing scores as values of a covariate have to be assigned to each category of the risk factor. Thus, the $p$-value of the trend test may depend on such scores. A collection of trend tests is formed by choosing various increasing scores. (Graubard and Korn [12]) A third example arises in genetic linkage and association studies. In linkage analysis to map quantitative trait loci using affected sib pairs, optimal tests are functions of the number of alleles shared identical-by-descent (IBD) by the two sibs. The IBD probabilities form a family of alternatives which are determined by genetic models. See, e.g., Whittemore and Tu [22] and Gastwirth and Freidlin [10]. In genetic association studies using case-parents trios, the optimal test depends on the mode of inheritance of the disease (recessive, dominant, or co-dominant disease). For complex diseases, the underlying genetic model is often not known. Using a single optimal test does not protect against a substantial loss of power under the worst situation, i.e., when a much different genetic model is the true one. (Zheng, Freidlin and Gastwirth [25])

Robust procedures have been developed and applied when the underlying model is unknown as discussed in Gastwirth [7, 8, 9], Birnbaum and Laska [2], Podgor, Gastwirth and Mehta [16], and Freidlin, Podgor and Gastwirth [5]. In this article, we review two useful robust tests. The first one is a linear combination of the two or three extreme optimal tests in a family of optimal statistics and the second one is a suitably chosen maximum statistic, i.e., the maximum of several of the optimum tests for specific models in the family. These two robust procedures are applied to genetic association using case-control studies and compared to other test statistics that are used in practice.

## 2. Robust procedures: A short review

Suppose we have a collection of alternative models $\{M_i, i \in I\}$ and the corresponding optimal (most powerful) test statistics $\{T_i : i \in I\}$ are obtained, where $I$ can be a finite set or an interval. Under the null hypothesis, assume that each of these test statistics is asymptotically normally distributed, i.e., $Z_i = [T_i - \mathrm{E}(T_i)]/\{\mathrm{Var}(T_i)\}^{1/2}$ converges in law to $N(0,1)$ where $\mathrm{E}(T_i)$ and $\mathrm{Var}(T_i)$ are the mean and the variance of $T_i$ under the null; suppose also that for any $i, j \in I$, $Z_i$ and $Z_j$ are jointly normal with the correlation $\rho_{ij}$. When $M_i$ is the true model, the optimal test $Z_i$ would be used. When the true model $M_i$ is unknown and the test $Z_j$ is used, assume the Pitman asymptotic relative efficiency (ARE) of $Z_j$ relative to $Z_i$ is $e(Z_j, Z_i) = \rho_{ij}^2$ for $i, j \in I$. These conditions are satisfied in many applications. (van Eeden [21] and Gross [13])

### 2.1. Maximin efficiency robust tests

When the true model is unknown and each model in the family is scientifically plausible, the minimum ARE compared to the optimum test for each model, $Z_i$, when $Z_j$ is used is given by $\inf_{i\in I} e(Z_j, Z_i)$ for $j \in I$. One robust test is to choose the optimal test $Z_l$ from the family $\{Z_i : i \in I\}$ which maximizes the minimum ARE, that is,

$$\inf_{i\in I} e(Z_l, Z_i) = \sup_{j\in I} \inf_{i\in I} e(Z_j, Z_i). \tag{2.1}$$

Under the null hypothesis, $Z_l$ converges in distribution to a standard normal random variable and, under the definition (2.1), is the most robust test in $\{Z_i : i \in I\}$. In practice, however, other tests have been studied which may have greater efficiency robustness.

Although a family of models are proposed based on scientific knowledge and the corresponding optimal tests can be obtained, all consistent tests with an asymptotically normal distribution can be used. Denote all these tests for the problem by $C$. The original family of test statistics can be expanded to $C$. The purpose is to find a test $Z$ from $C$, rather than from the original family $\{Z_i : i \in I\}$, such that

$$\inf_{i\in I} e(Z, Z_i) = \sup_{Z\in C} \inf_{i\in I} e(Z, Z_i). \tag{2.2}$$

The test $Z$ satisfying (2.2) is called maximin efficiency robust test (MERT). (Gastwirth [7]) When the family $C$ is restricted to the convex linear combinations of $\{Z_i : i \in I\}$, the resulting robust test is denoted as $Z_{\text{MERT}}$. Since $\{Z_i : i \in I\} \subset C$,

$$\sup_{Z\in C} \inf_{i\in i} e(Z, Z_i) \geq \sup_{j\in I} \inf_{i\in I} e(Z_j, Z_i).$$

Assuming that $\inf_{i,j\in I} \rho_{ij} \geq \epsilon > 0$, Gastwirth [7] proved that $Z_{\text{MERT}}$ uniquely exists and can be written as a closed convex combination of optimal tests $Z_i$ in the family $\{Z_i : i \in i\}$. Although a simple algorithm when $C$ is the class of linear combination of $\{Z_i : i \in I\}$ was given in Gastwirth [9] (see also Zucker and Lakatos [27]), the computation of $Z_{\text{MERT}}$ is more complicated as it is related to quadratic programming algorithms. (Rosen [18]) For many applications, $Z_{\text{MERT}}$ can be easily written as a linear convex combination of two or three optimal tests in $\{Z_i : i \in I\}$ including the extreme pair defined as follows: two optimal tests $Z_s, Z_t \in \{Z_i : i \in I\}$ are called extreme pair if $\rho_{st} = \text{corr}_{H_0}(Z_s, Z_t) = \inf_{i,j\in I} \rho_{ij} > 0$. Define a new test statistic $Z_{st}$ based on the extreme pair as

$$Z_{st} = \frac{Z_s + Z_t}{[2(1+\rho_{st})]^{1/2}}, \tag{2.3}$$

which is the MERT for the extreme pair. A necessary and sufficient condition for $Z_{st}$ to be $Z_{\text{MERT}}$ for the whole family $\{Z_i : i \in I\}$ is given, see Gastwirth [8], by

$$\rho_{si} + \rho_{it} \geq 1 + \rho_{st}, \quad \text{for all } i \in I. \tag{2.4}$$

Under the null hypothesis, $Z_{\text{MERT}}$ is asymptotically $N(0,1)$. The ARE of the MERT given by (2.4) is $(1+\rho_{st})/2$. To find the MERT, the null correlations $\rho_{ij}$ need to be obtained and the pair is the extreme pair for which $\rho_{ij}$ is smallest.

### 2.2. Maximum tests

The robust test $Z_{\mathrm{MERT}}$ is a linear combination of the optimal test statistics and with modern computers it is useful to extend the family $C$ of possible tests to include non-linear functions of the $Z_i$. A natural non-linear robust statistic is the maximum over the extreme pair $(Z_s, Z_t)$ or the triple $(Z_s, Z_u, Z_t)$ for the entire family (Freidlin et al. [5]), i.e.,

$$Z_{\mathrm{MAX2}} = \max(Z_s, Z_t) \;\text{ or }\; Z_{\mathrm{MAX3}} = \max(Z_s, Z_u, Z_t).$$

There are several choices for $Z_u$ in $Z_{\mathrm{MAX3}}$, e.g., $Z_u = Z_{st}$ (MERT for the extreme pair or entire family). As when obtaining the MERT, the correlation matrix $\{\rho_{ij}\}$ guides the choice of $Z_u$ to be used in MAX3, e.g., it has equal correlation with the extreme tests. A more complicated maximum test statistic is to take the maximum over the entire family $Z_{\mathrm{MAX}} = \max_{i\in I} Z_i$ or $Z_{\mathrm{MAX}} = \max_{i\in C} Z_i$. $Z_{\mathrm{MAX}}$ was considered by Davies [4] for some non-standard hypothesis testing whose critical value has to be determined by approximation of its upper bound. In a recent study of several applications in genetic association and linkage analysis, Zheng and Chen [24] showed that $Z_{\mathrm{MAX3}}$ and $Z_{\mathrm{MAX}}$ have similar power performance in these applications. Moreover, $Z_{\mathrm{MAX2}}$ or $Z_{\mathrm{MAX3}}$ are much easier to compute than $Z_{\mathrm{MAX}}$. Hence, in the next section, we only consider the two maximum tests $Z_{\mathrm{MAX2}}$ and $Z_{\mathrm{MAX3}}$. The critical values for the maximum test statistics can be found by simulation under the null hypothesis as any two or three optimal statistics in $\{Z_i : i \in I\}$ follow multivariate normal distributions with correlation matrix $\{\rho_{ij}\}$. For example, given the data, $\rho_{st}$ can be calculated. Generating a bivariate normal random variable $(Z_{sj}, Z_{tj})$ with the correlation matrix $\{p_{st}\}$ for $j = 1, \ldots, B$. For each $j$, $Z_{MAX2}$ is obtained. Then an empirical distribution for $Z_{\mathrm{MAX2}}$ can be obtained using these $B$ simulated maximum statistics, from which we can find the critical values. In some applications, if the null hypothesis does not depend on any nuisance parameters, the distribution of $Z_{\mathrm{MAX2}}$ or $Z_{\mathrm{MAX3}}$ can be simulated exactly without the correlation matrix, e.g., Zheng and Chen [24].

### 2.3. Comparison of MERT and MAX

Usually, $Z_{\mathrm{MERT}}$ is easier to compute and use than $Z_{\mathrm{MAX2}}$ (or $Z_{\mathrm{MAX3}}$). Intuitively, however, $Z_{\mathrm{MAX3}}$ should have greater efficiency robustness than $Z_{\mathrm{MERT}}$ when the range of models is wide. The selection of the robust test depends on the minimum correlation $\rho_{st}$ of the entire family of optimal tests. Results from Freidlin et al. [5] showed that when $\rho_{st} \geq 0.75$, MERT and MAX2 (MAX3) have similar power; thus, the simpler MERT can be used. For example, when $\rho_{st} = 0.75$, the ARE of MERT relative to the optimal test for any model in the family is at least 0.875. When $\rho_{st} < 0.50$, MAX2 (MAX3) is noticeably more powerful than the simple MERT. Hence, MAX2 (MAX3) is recommended. For example, in genetic linkage analysis using affected sib pairs, the minimum correlation is greater than 0.8, and the MERT, MAX2, and MAX3 have similar power. (Whittemore and Tu [22] and Gastwirth and Freidlin [10]) For analysis of case-parents data in genetic association studies where the mode of inheritance can range from pure recessive to pure dominant, the minimum correlation is less than 0.33, and then the MAX3 has desirable power robustness for this problem. (Zheng et al. [25])

## 3. Genetic association using case-control studies

### *3.1. Background*

It is well known that association studies testing linkage disequilibrium are more powerful than linkage analysis to detect small genetic effects on traits. (Risch and Merikangas [17]) Moreover association studies using cases and controls are easier to conduct as parental genotypes are not required.

Assume that cases are sampled from the study population and that controls are independently sampled from the general population without disease. Cases and controls are not matched. Each individual is genotyped with one of three genotypes $MM$, $MN$ and $NN$ for a marker with two alleles $M$ and $N$. The data obtained in case-control studies can be displayed as in Table 1 (genotype-based) or as in Table 2 (allele-based).

Define three penetrances as $f_0 = \Pr(\text{case}|NN)$, $f_1 = \Pr(\text{case}|NM)$, and $f_2 = \Pr(\text{case}|MM)$, which are the disease probabilities given different genotypes. The prevalence of disease is denoted as $D = \Pr(\text{case})$. The probabilities for genotypes $(NN, NM, MM)$ in cases and controls are denoted by $(p_0, p_1, p_2)$ and $(q_0, q_1, q_2)$, respectively. The probabilities for genotypes $(NN, NM, MM)$ in the general population are denoted as $(g_0, g_1, g_2)$. The following relations can be obtained.

$$p_i = \frac{f_i g_i}{D} \text{ and } q_i = \frac{(1-f_i)g_i}{1-D} \text{ for } i = 0, 1, 2. \tag{3.1}$$

Note that, in Table 1, $(r_0, r_1, r_2)$ and $(s_0, s_1, s_2)$ follow multinomial distributions $mul(R; p_0, p_1, p_2)$ and $mul(S; q_0, q_1, q_2)$, respectively. Under the null hypothesis of no association between the disease and the marker, $p_i = q_i = g_i$ for $i = 0, 1, 2$. Hence, from (3.1), the null hypothesis for Table 1 is equivalent to $H_0 : f_0 = f_1 = f_2 = D$. Under the alternative, penetrances are different as one of two alleles is a risk allele, say, $M$. In genetic analysis, three genetic models (mode of inheritance) are often used. A model is recessive (rec) when $f_0 = f_1$, additive (add) when $f_1 = (f_0 + f_2)/2$, and dominant (dom) when $f_1 = f_2$. For recessive and dominant models, the number of columns in Table 1 can be reduced. Indeed, the columns with $NN$ and $NM$ ($NM$ and $MM$) can be collapsed for recessive (dominant) model. Testing association using Table 2 is simpler but Sasieni [19] showed that genotype based analysis is preferable unless cases and controls are in Hardy–Weinberg Equilibrium.

Table 1
*Genotype distribution for case-control studies*

| | $NN$ | $NM$ | $MM$ | Total |
|---|---|---|---|---|
| Case | $r_0$ | $r_1$ | $r_2$ | $r$ |
| Control | $s_0$ | $s_1$ | $s_2$ | $s$ |
| Total | $n_0$ | $n_1$ | $n_2$ | $n$ |

Table 2
*Allele distribution for case-control studies*

| | $N$ | $M$ | Total |
|---|---|---|---|
| Case | $2r_0 + r_1$ | $r_1 + 2r_2$ | $2r$ |
| Control | $2s_0 + s_1$ | $s_1 + 2s_2$ | $2s$ |
| Total | $2n_0 + n_1$ | $n_1 + 2n_2$ | $2n$ |

### 3.2. Test statistics

For the $2 \times 3$ table (Table 1), a chi-squared test with 2 degrees of freedom (df) can be used. (Gibson and Muse [11]) This test is independent of the underlying genetic model. Note that, under the alternative when $M$ is the risk allele, the penetrances have a natural order: $f_0 \le f_1 \le f_2$ (at least one inequality hold). The Cochran-Armitage (CA) trend test (Cochran [3] and Armitage [1]) taking into account the natural order should be more powerful than the chi-squared test as the trend test has 1 df.

The CA trend test can be obtained as a score test under the logistic regression model with genotype as a covariate, which is coded using scores $\mathbf{x} = (x_0, x_1, x_1)$ for $(NN, NM, MM)$, where $x_0 \le x_1 \le x_2$. The trend test can be written as (Sasieni [19])

$$Z_{\mathbf{x}} = \frac{n^{1/2} \sum_{i=0}^{2} x_i(sr_i - rs_i)}{\{rs[n \sum_{i=0}^{2} x_i^2 n_i - (\sum_{i=0}^{2} x_i n_i)^2]\}^{1/2}}.$$

Since the trend test is invariant to linear transformations of $\mathbf{x}$, without loss of generality, we use the scores $\mathbf{x} = (0, x, 1)$ with $0 \le x \le 1$ and denote $Z_{\mathbf{x}}$ as $Z_x$. Under the null hypothesis, $Z_x$ has an asymptotic normal distribution $N(0,1)$. When $M$ is a risk allele, a one-sided test is used. Otherwise, a two-sided test should be used. Results from Sasieni [19] and Zheng, Freidlin, Li and Gastwirth [26] showed that the optimal choices of $x$ for recessive, additive and dominant models are $x = 0$, $x = 1/2$, and $x = 1$, respectively. That is, $Z_0$, $Z_{1/2}$ or $Z_1$ is an asymptotically most powerful test when the genetic model is recessive, additive or dominant. The tests using other values of $x$ are optimal for penetrances in the range $0 < f_0 \le f_1 \le f_2 < 1$.

For complex diseases, the genetic model is not known a priori. The optimal test $Z_x$ cannot be used directly as a substantial loss of power may occur when $x$ is misspecified. Applying the robust procedures introduced in Section 2, we have three genetic models and the collection of all consistent tests $C = \{Z_x : x \in [0,1]\}$. To find a robust test, we need to evaluate the null correlations. Denote these as $\text{corr}_{H_0}(Z_{x_1}, Z_{x_2}) = \rho_{x_1,x_2}$. From appendix C of Freidlin, Zheng, Li and Gastwirth [6],

$$\begin{aligned}
\rho_{0,1/2} &= \frac{p_0(p_1 + 2p_2)}{\{p_0(1-p_0)\}^{1/2}\{(p_1+2p_2)p_0 + (p_1+2p_0)p_2\}^{1/2}}, \\
\rho_{0,1} &= \frac{p_0 p_2}{\{p_0(1-p_0)\}^{1/2}\{p_2(1-p_2)\}^{1/2}}, \\
\rho_{1/2,1} &= \frac{p_2(p_1 + 2p_0)}{\{p_2(1-p_2)\}^{1/2}\{(p_1+2p_2)p_0 + (p_1+2p_0)p_2\}^{1/2}}.
\end{aligned}$$

Although the null correlations are functions of the unknown parameters $p_i, i = 0, 1, 2$, it can be shown analytically that $\rho_{0,1} < \rho_{0,1/2}$ and $\rho_{0,1} < \rho_{1/2,1}$. Note that if the above analytical results were not available, the $p_i$ would be estimated by substituting the observed data $\hat{p}_i = n_i/n$ for $p_i$. Here the minimum correlation among the three optimal tests occurs when $Z_0$ and $Z_1$ is the extreme pair for the three genetic models. Freidlin et al. [6] also proved analytically that the condition (2.4) holds. Hence, $Z_{\text{MERT}} = (Z_0 + Z_1)/\{2(1 + \hat{\rho}_{0,1})\}^{1/2}$ is the MERT for the whole family $C$, where $\hat{\rho}_{0,1}$ is obtained when the $p_i$ are replaced by $n_i/n$. The two maximum tests can be written as $Z_{\text{MAX2}} = \max(Z_0, Z_1)$ and $Z_{\text{MAX3}} = \max(Z_0, Z_{1/2}, Z_1)$. When the risk allele is unknown, $Z_{\text{MAX2}} = \max(|Z_0|, |Z_1|)$ and $Z_{\text{MAX3}} = \max(|Z_0|, |Z_{1/2}|, |Z_1|)$. Although we considered three genetic models, the family of genetic models for case-control studies can be extended by defining a genetic model as penetrances restricted to the family $\{(f_0, f_1, f_2) : f_0 \leq f_1 \leq f_2\}$. Three genetic models are contained in this family as the two boundaries and one middle ray of this family. The statistics $Z_{\text{MERT}}$ and $Z_{\text{MAX3}}$ are also the corresponding robust statistics for this larger family (see, e.g., Freidlin et al. [6] and Zheng et al. [26]).

In analysis of case-control data for genetic association, two other tests are also currently used. However, their robustness and efficiency properties have not been compared to MERT and MAX. The first one is the chi-squared test for the $2 \times 3$ contingency table (Table 1), denoted as $\chi_2^2$. (Gibson and Muse [11]) Under the null hypothesis, it has a chi-squared distribution with 2 df. The second test, denoted as $Z_{\text{P}}$, is based on the product of two different tests: (a) the allele association (AA) test and (b) the Hardy-Weinberg disequilibrium (HWD) test. (Hoh, Wile and Ott [15] and Song and Elston [20]) The AA test is a chi-squares test for the $2 \times 2$ table given in Table 2, which is written as

$$\chi_{\text{AA}}^2 = \frac{2n[(2r_0 + r_1)(s_1 + 2s_2) - (2s_0 + s_1)(r_1 + 2r_2)]^2}{4rs(2n_0 + n_1)(n_1 + 2n_2)}.$$

The HWD test detects the deviation from Hardy–Weinberg equilibrium (HWE) in cases. Assume the allele frequency of $M$ is $p = \Pr(M)$. Using cases, the estimation of $p$ is $\hat{p} = (r_1 + 2r_2)/(2r)$. Let $\hat{q} = 1 - \hat{p}$ be the estimation of allele frequency for $N$. Under the null hypothesis of HWE, the expected number of genotypes can be written as $\mathrm{E}(NN) = r\hat{q}^2$, $\mathrm{E}(NM) = r2\hat{p}\hat{q}$ and $\mathrm{E}(MM) = r\hat{p}^2$, respectively. Hence, a chi-squared test for HWE is

$$\chi_{\text{HWD}}^2 = \frac{(r_0 - \mathrm{E}(NN))^2}{\mathrm{E}(NN)} + \frac{(r_1 - \mathrm{E}(NM))^2}{\mathrm{E}(NM)} + \frac{(r_2 - \mathrm{E}(MM))^2}{\mathrm{E}(MM)}.$$

The product test, proposed by Hoh et al. [15], is $T_{\text{P}} = \chi_{\text{AA}}^2 \times \chi_{\text{HWD}}^2$. They noticed that the power performances of these two statistics are complementary. Thus, the product should retain reasonable power as one of the tests has high power when the other does not. Consequently, for a comprehensive comparison, we also consider the maximum of them, $T_{\text{MAX}} = \max(\chi_{\text{AA}}^2, \chi_{\text{HWD}}^2)$. Given the data, the critical values of $T_{\text{P}}$ and $T_{\text{MAX}}$ can be obtained by a permutation procedure as their asymptotic distributions are not available. (Hoh et al. [15]) Note that $T_{\text{P}}$ was originally proposed by Hoh et al. [15] as a test statistic for multiple gene selection and was modified by Song and Elston [20] for use as a test statistic for a single gene.

### *3.3. Power comparison*

We conducted a simulation study to compare the power performance of the test statistics. The test statistics were (a) the optimal trend tests for the three genetic models, $Z_0$, $Z_{1/2}$ and $Z_1$, (b) MERT $Z_{\text{MERT}}$, (c) maximum tests $Z_{\text{MAX2}}$ and $Z_{\text{MAX3}}$, (d) the product test $T_{\text{P}}$, (e) $T_{\text{MAX}}$, and (f) $\chi_2^2$.

In the simulation a two-sided test was used. We assumed that the allele frequency $p$ and the baseline penetrance $f_0$ are known ($f_0 = .01$). Note that, in practice, the allele frequency and penetrances are unknown. However, they can be estimated empirically (e.g., Song and Elston [20] and Wittke-Thompson, Pluzhnikov and Cox [23]). In our simulation the critical values for all test statistics are simulated under the null hypothesis. Thus, we avoid using asymptotic distributions for the test statistics. The Type I errors for all tests are expected to be close to the nominal level $\alpha = 0.05$ and the powers of all tests are comparable. When HWE holds, the probabilities $(p_0, p_1, p_2)$ for cases and $(q_0, q_1, q_2)$ for controls can be calculated using (3.1) under the null and alternative hypotheses, where $(g_0, g_1, g_2) = (q^2, 2pq, q^2)$ and $(f_1, f_2)$ are specified by the null or alternative hypotheses and $D = \sum f_i g_i$. After calculating $(p_0, p_1, p_2)$ and $(q_0, q_1, q_2)$ under the null hypothesis, we first simulated the genotype distributions $(r_0, r_1, r_2) \sim mul(R; p_0, p_1, p_2)$ and $(s_0, s_1, s_2) \sim mul(S; q_0, q_1, q_2)$ for cases and controls, respectively (see Table 1). When HWE does not hold, we assumed a mixture of two populations with two different allele frequencies $p_1$ and $p_2$. Hence, we simulated two independent samples with different allele frequencies for cases (and controls) and combined these two samples for cases (and for controls). Thus, cases (controls) contain samples from a mixture of two populations with different allele frequencies. When $p$ is small, some counts can be zero. Therefore, we added 1/2 to the count of each genotype in cases and controls in all simulations.

To obtain the critical values, a simulation under the null hypothesis was done with 200,000 replicates. For each replicate, we calculated the test statistics. For each test statistic, we used its empirical distribution function based on 200,000 replicates to calculate the critical value for $\alpha = 0.05$. The alternatives were chosen so that the power of the optimal test $Z_0, Z_{1/2}, Z_1$ was near 80% for the recessive, additive, dominant models, respectively. To determine the empirical power, 10,000 replicates were simulated using multinomial distributions with the above probabilities.

To calculate $Z_{\text{MERT}}$, the correlation $\rho_{0,1}$ was estimated using the simulated data. In Table 3, we present the mean of correlation matrix using 10,000 replicates when $r = s = 250$. The three correlations $\rho_{0,1/2}, \rho_{0,1}, \rho_{1/2,1}$ were estimated by replacing $p_i$ with $n_i/n, i = 0, 1, 2$ using the data simulated under the null and alternatives and various models. The null and alternative hypotheses used in Table 3 were also used to simulate critical values and powers (Table 4). Note that the minimum correlation $\rho_{0,1}$ is less than .50. Hence, the $Z_{\text{MAX3}}$ should have greater efficiency robustness than $Z_{\text{MERT}}$. However, when the dominant model can be eliminated based on prior scientific knowledge (e.g. the disease often skips generations), the correlation between $Z_0$ and $Z_{1/2}$ optimal for the recessive and additive models would be greater than .75. Thus, for these two models, $Z_{\text{MERT}}$ should have comparable power to $Z_{\text{MAX2}} = \max(|Z_0|, |Z_{1/2}|)$ and is easier to use.

The correlation matrices used in Table 4 for $r \neq s$ and Table 5 for mixed samples are not presented as they did not differ very much from those given in Table 3. Tables 4 and 5 present simulation results where all three genetic models are plausible. When HWE holds Table 4 shows that the Type I error is indeed close to the $\alpha = 0.05$ level. Since the model is not known, the minimum power across three genetic models is written in bold type. A test with the maximum of the minimum power among all test statistic has the most power robustness. Our comparison fo-

TABLE 3
*The mean correlation matrices of three optimal test statistics based on 10,000 replicates when HWE holds*

| | $p$ | | | | | | | | |
|---|---|---|---|---|---|---|---|---|---|
| | .1 | | | .3 | | | .5 | | |
| Model | $\rho_{0,1/2}$ | $\rho_{0,1}$ | $\rho_{1/2,1}$ | $\rho_{0,1/2}$ | $\rho_{0,1}$ | $\rho_{1/2,1}$ | $\rho_{0,1/2}$ | $\rho_{0,1}$ | $\rho_{1/2,1}$ |
| null | .97 | .22 | .45 | .91 | .31 | .68 | .82 | .33 | .82 |
| rec | .95 | .34 | .63 | .89 | .36 | .74 | .81 | .37 | .84 |
| add | .96 | .23 | .48 | .89 | .32 | .71 | .79 | .33 | .84 |
| dom | .97 | .21 | .44 | .90 | .29 | .69 | .79 | .30 | .83 |

The same models (rec,add,dom) are used in Table 4 when $r = s = 250$.

Table 4
*Power comparison when HWE holds in cases and controls under three genetic models with $\alpha = .05$*

| | | Test statistics | | | | | | | | |
|---|---|---|---|---|---|---|---|---|---|---|
| $p$ | Model | $Z_0$ | $Z_{1/2}$ | $Z_1$ | $Z_{\mathrm{MERT}}$ | $Z_{\mathrm{MAX2}}$ | $Z_{\mathrm{MAX3}}$ | $T_{\mathrm{MAX}}$ | $T_{\mathrm{P}}$ | $\chi_2^2$ |
| | | | | | $r = 250, s = 250$ | | | | | |
| .1 | null | .058 | .048 | .050 | .048 | .047 | .049 | .049 | .049 | .053 |
| | rec | .813 | **.364** | **.138** | **.606** | **.725** | **.732** | .941 | .862 | **.692** |
| | add | .223 | .813 | .802 | .733 | .782 | .800 | .705 | **.424** | .752 |
| | dom | **.108** | .796 | .813 | .635 | .786 | .795 | **.676** | .556 | .763 |
| .3 | null | .051 | .052 | .051 | .054 | .052 | .052 | .053 | .049 | .050 |
| | rec | .793 | **.537** | **.178** | .623 | **.714** | **.726** | .833 | .846 | **.691** |
| | add | .433 | .812 | .768 | .786 | .742 | .773 | .735 | **.447** | .733 |
| | dom | **.133** | .717 | .809 | **.621** | .737 | .746 | **.722** | .750 | .719 |
| .5 | null | .049 | .047 | .051 | .047 | .047 | .047 | .050 | .051 | .050 |
| | rec | .810 | .662 | **.177** | .644 | .738 | .738 | .772 | .813 | .709 |
| | add | .575 | .802 | .684 | .807 | .729 | .760 | **.719** | **.450** | .714 |
| | dom | **.131** | **.574** | .787 | **.597** | **.714** | **.713** | .747 | .802 | **.698** |
| | | | | | $r = 50, s = 250$ | | | | | |
| .1 | null | .035 | .052 | .049 | .052 | .051 | .053 | .045 | .051 | .052 |
| | rec | .826 | **.553** | **.230** | .757 | .797 | .802 | .779 | .823 | .803 |
| | add | .250 | .859 | .842 | .773 | .784 | .795 | .718 | **.423** | .789 |
| | dom | **.114** | .807 | .814 | **.658** | **.730** | **.734** | **.636** | .447 | **.727** |
| .3 | null | .048 | .048 | .048 | .048 | .050 | .050 | .050 | .048 | .049 |
| | rec | .836 | **.616** | **.190** | .715 | .787 | .787 | .733 | .752 | .771 |
| | add | .507 | .844 | .794 | .821 | .786 | .813 | .789 | **.500** | .778 |
| | dom | **.171** | .728 | .812 | **.633** | **.746** | **.749** | **.696** | .684 | **.729** |
| .5 | null | .049 | .046 | .046 | .046 | .047 | .047 | .046 | .049 | .046 |
| | rec | .838 | .692 | **.150** | .682 | .771 | .765 | .705 | .676 | .748 |
| | add | .615 | .818 | .697 | .820 | .744 | .780 | .743 | **.493** | .728 |
| | dom | **.151** | **.556** | .799 | **.565** | **.710** | **.708** | **.662** | .746 | **.684** |

cuses on the test statistics: $Z_{\mathrm{MAX3}}$, $T_{\mathrm{MAX}}$, $T_{\mathrm{P}}$ and $\chi_2^2$. Our results show that $T_{\mathrm{MAX}}$ has greater efficiency robustness than $T_{\mathrm{P}}$ while $Z_{\mathrm{MAX3}}$, $T_{\mathrm{MAX}}$ and $\chi_2^2$ have similar minimum powers. Notice that $Z_{\mathrm{MAX3}}$ is preferable to $\chi_2^2$ although the difference in minimum powers depends on the allele frequency. When HWE does not hold, $Z_{\mathrm{MAX3}}$ still possesses its efficiency robustness, but $T_{\mathrm{MAX}}$ and $T_{\mathrm{P}}$ do not perform as well. Thus, population stratification affects their performance. $Z_{\mathrm{MAX3}}$ also remains more robust than $\chi_2^2$ even when HWE does not hold. From both Table 4 and Table 5, $\chi_2^2$ is more powerful than $Z_{\mathrm{MERT}}$ except for the additive model. However, when the genetic model is known, the corresponding optimal CA trend test is more powerful than $\chi_2^2$ with 2 df.

From Tables 4 and 5, one sees that the robust test $Z_{\mathrm{MAX3}}$ tends to be more powerful than $\chi_2^2$ under the various scenarios we simulated. Further comparisons of these two test statistics using $p$-values and the same data given in Table 4 (when $r = s = 250$) are reported in Table 6. Following Zheng et al. [25], which reported a matched study, where both tests are applied to the same data, the $p$-values for each test are grouped as $< .01$, $(.01, .05)$, $(.05, .10)$, and $> .10$. Cross classification of the $p$-values are given in Table 6 for allele frequencies $p = .1$ and .3 under all three genetic models. Table 6 is consistent with the results of Tables 4 and 5, i.e., $Z_{\mathrm{MAX3}}$ is more powerful than the chi-squared test with 2 degrees of freedom when the genetic model is unknown. This is seen by comparing the counts at the upper right corner with the counts at lower left corner. When the counts at the upper right corner are greater than the corresponding counts at the lower left corner, $Z_{\mathrm{MAX3}}$ usually has smaller $p$-values than $\chi_2^2$. In particular, we compare two tests with $p$-values $< .01$ versus $p$-values in $(.01, .05)$. Notice that in most situations

Table 5

*Power comparison when HWE does not hold in cases and controls under three genetic models (Mixed samples with different allele frequencies $(p_1, p_2)$ and sample sizes $(R_1, S_1)$ and $(R_2, S_2)$ with $r = R_1 + R_2$, $s = S_1 + S_2$, and $\alpha = .05$).*

| | | Test statistics | | | | | | | | |
|---|---|---|---|---|---|---|---|---|---|---|
| $(p_1, p_2)$ | Model | $Z_0$ | $Z_{1/2}$ | $Z_1$ | $Z_{\mathrm{MERT}}$ | $Z_{\mathrm{MAX2}}$ | $Z_{\mathrm{MAX3}}$ | $T_{\mathrm{MAX}}$ | $T_{\mathrm{P}}$ | $\chi_2^2$ |
| | | $R_1 = 250, S_1 = 250$ and $R_2 = 100, S_2 = 100$ | | | | | | | | |
| (.1,.4) | null | .047 | .049 | .046 | .050 | .047 | .046 | .048 | .046 | .048 |
| | rec | .805 | **.519** | **.135** | .641 | .747 | .744 | .936 | .868 | .723 |
| | add | .361 | .817 | .771 | .776 | .737 | .757 | .122 | .490 | .688 |
| | dom | **.098** | .715 | .805 | **.586** | **.723** | **.724** | **.049** | **.097** | **.681** |
| (.1,.5) | null | .047 | .050 | .046 | .048 | .046 | .046 | .052 | .050 | .052 |
| | rec | .797 | **.537** | **.121** | .620 | .746 | .750 | .920 | .839 | .724 |
| | add | .383 | .794 | .732 | .764 | **.681** | .709 | .033 | .620 | **.631** |
| | dom | **.095** | .695 | .812 | **.581** | .697 | **.703** | **.001** | **.133** | .670 |
| (.2,.5) | null | .048 | .052 | .052 | .054 | .052 | .052 | .050 | .047 | .050 |
| | rec | .816 | **.576** | **.157** | .647 | .760 | .749 | .889 | .881 | .725 |
| | add | .417 | .802 | .754 | .782 | **.726** | .743 | .265 | .365 | **.684** |
| | dom | **.112** | .679 | .812 | **.600** | .729 | **.715** | **.137** | **.122** | .696 |
| | | $R_1 = 30, S_1 = 150$ and $R_2 = 20, S_2 = 100$ | | | | | | | | |
| (.1,.4) | null | .046 | .048 | .048 | .046 | .048 | .046 | .055 | .045 | .048 |
| | rec | .847 | **.603** | **.163** | .720 | .807 | .799 | .768 | .843 | .779 |
| | add | .387 | .798 | .762 | .749 | .725 | .746 | .449 | .309 | **.691** |
| | dom | **.139** | .733 | .810 | **.603** | **.721** | **.728** | **.345** | **.223** | .699 |
| (.1,.5) | null | .053 | .055 | .050 | .053 | .053 | .053 | .059 | .047 | .054 |
| | rec | .816 | **.603** | **.139** | .688 | .776 | .780 | .697 | .827 | .741 |
| | add | .415 | .839 | .797 | .798 | .763 | .781 | .276 | .358 | **.703** |
| | dom | **.120** | .726 | .845 | **.612** | **.750** | **.752** | **.149** | **.133** | .716 |
| (.2,.5) | null | .047 | .048 | .050 | .051 | .048 | .048 | .049 | .050 | .044 |
| | rec | .858 | **.647** | **.167** | .722 | .808 | .804 | .768 | .852 | .783 |
| | add | .472 | .839 | .790 | .811 | .776 | .799 | .625 | **.325** | .740 |
| | dom | **.139** | .708 | .815 | **.614** | **.741** | **.743** | **.442** | .390 | **.708** |

the number of times $\chi_2^2$ has a $p$-value $< .01$ and $Z_{\mathrm{MAX3}}$ has a $p$-value in $(.01, .05)$ is much less than the corresponding number of times when $Z_{\mathrm{MAX3}}$ has a $p$-value $< .01$ and $\chi_2^2$ has a $p$-value in $(.01, .05)$. For example, when $p = .3$ and the additive model holds, there are 289 simulated datasets where $\chi_2^2$ has a $p$-value in (.01,.05) while $Z_{\mathrm{MAX3}}$ has a $p$-value $< .01$ versus only 14 such datasets when $Z_{\mathrm{MAX3}}$ has a $p$-value in (.01,.05) while $\chi_2^2$ has a $p$-value $< .01$. The only exception occurs at the recessive model under which they have similar counts (165 vs. 140). Combining results from Tables 4 and 6, $Z_{\mathrm{MAX3}}$ is more powerful than $\chi_2^2$, but the difference of power between $Z_{\mathrm{MAX3}}$ and $\chi_2^2$ is usually less than 5% in the simulation. Hence $\chi_2^2$ is also an efficiency robust test, which is very useful for genome-wide association studies, where hundreds of thousands of tests are performed.

From prior studies of family pedigrees one may know whether the disease skips generations or not. If it does, the disease is less likely to follow a pure-dominant model. Thus, when genetic evidence strongly suggests that the underlying genetic model is between the recessive and additive inclusive, we compared the performance of tests $Z_{\mathrm{MERT}} = (Z_0 + Z_{1/2})/\{2(1 + \hat{\rho}_{0,1/2})\}^{1/2}$, $Z_{\mathrm{MAX2}} = \max(Z_0, Z_{1/2})$, and $\chi_2^2$. The results are presented in Table 7. The alternatives used in Table 4 for rec and add with $r = s = 250$ were also used to obtain Table 7. For a family with the recessive and additive models, the minimum correlation is increased compared to the family with three genetic models (rec, add and dom). For example, from Table 3, the minimum correlation with the family of three models that ranges from .21 to .37 is increased to the range of .79 to .97 with only two models. From Table 7, under the recessive and additive models, while $Z_{\mathrm{MAX2}}$ remains more powerful than other[1] is dominant and the other two are semi-dominant, all with $f_2 = .019$.
other [1]: $f_1 = .019$. other[2]: $f_1 = .017$. other[3]: $f_1 = .015$.

TABLE 6
*Matched p-value comparison of $Z_{MAX3}$ and $\chi^2_2$ when HWE holds in cases and controls under three genetic models (Sample sizes $r = s = 250$ and 5,000 replications)*

| $p$ | Model | $Z_{MAX3}$ | $\chi^2_2$ $< .01$ | $.01 - .05$ | $.05 - .10$ | $> .10$ |
|---|---|---|---|---|---|---|
| .10 | rec | $< .01$ | 2069 | 165 | 0 | 0 |
| | | $.01 - .05$ | 140 | 1008 | 251 | 0 |
| | | $.05 - .10$ | 7 | 52 | 227 | 203 |
| | | $> .10$ | 1 | 25 | 47 | 805 |
| | add | $< .01$ | 2658 | 295 | 0 | 0 |
| | | $.01 - .05$ | 44 | 776 | 212 | 0 |
| | | $.05 - .10$ | 3 | 23 | 159 | 198 |
| | | $> .10$ | 0 | 5 | 16 | 611 |
| | dom | $< .01$ | 2785 | 214 | 0 | 0 |
| | | $.01 - .05$ | 80 | 712 | 214 | 0 |
| | | $.05 - .10$ | 10 | 42 | 130 | 169 |
| | | $> .10$ | 1 | 14 | 27 | 602 |
| .30 | rec | $< .01$ | 2159 | 220 | 0 | 0 |
| | | $.01 - .05$ | 85 | 880 | 211 | 0 |
| | | $.05 - .10$ | 6 | 44 | 260 | 157 |
| | | $> .10$ | 2 | 33 | 40 | 903 |
| | add | $< .01$ | 2485 | 289 | 0 | 0 |
| | | $.01 - .05$ | 14 | 849 | 229 | 0 |
| | | $.05 - .10$ | 1 | 8 | 212 | 215 |
| | | $> .10$ | 0 | 1 | 6 | 691 |
| | dom | $< .01$ | 2291 | 226 | 0 | 0 |
| | | $.01 - .05$ | 90 | 894 | 204 | 0 |
| | | $.05 - .10$ | 7 | 52 | 235 | 160 |
| | | $> .10$ | 0 | 26 | 49 | 766 |

TABLE 7
*Power comparison when HWE holds in cases and controls assuming two genetic models (rec and add) based on 10,000 replicates ($r = s = 250$ and $f_0 = .01$)*

| | $p$ | | | | | | | | |
|---|---|---|---|---|---|---|---|---|---|
| | .1 | | | .3 | | | .5 | | |
| Model | $Z_{MERT}$ | $Z_{MAX2}$ | $\chi^2_2$ | $Z_{MERT}$ | $Z_{MAX2}$ | $\chi^2_2$ | $Z_{MERT}$ | $Z_{MAX2}$ | $\chi^2_2$ |
| null | .046 | .042 | .048 | .052 | .052 | .052 | .053 | .053 | .053 |
| rec | .738 | .729 | .703 | .743 | .732 | .687 | .768 | .766 | .702 |
| add | .681 | .778 | .778 | .714 | .755 | .726 | .752 | .764 | .728 |
| other[1] | .617 | .830 | .836 | .637 | .844 | .901 | .378 | .547 | .734 |
| other[2] | .513 | .675 | .677 | .632 | .774 | .803 | .489 | .581 | .652 |
| other[3] | .385 | .465 | .455 | .616 | .672 | .655 | .617 | .641 | .606 |

$\chi^2_2$ and $Z_{MERT}$, the difference in minimum power is much less than in the previous simulation study (Table 4). Indeed, when studying complex common diseases where the allele frequency is thought to be fairly high, $Z_{MAX2}$ and $Z_{MERT}$ have similar power. Thus, when a genetic model is between the recessive and additive models inclusive, MAX2 and MERT should be used. In Table 7, some other models were also included in simulations when we do not have sound genetic knowledge to eliminate the dominant model. In this case, MAX2 and MERT lose some efficiency compared to $\chi^2_2$. However, MAX3 still has greater efficiency robustness than other tests. In particular, MAX3 is more powerful than $\chi^2_2$ (not reported) as in Table 4. Thus, MAX3 should be used when prior genetic studies do not justify excluding one of the basic three models.

## 4. Discussion

In this article, we review robust procedures for testing hypothesis when the underlying model is unknown. The implementation of these robust procedures is illustrated by applying them to testing genetic association in case-control studies. Simulation studies demonstrated the usefulness of these robust procedures when the underlying genetic model is unknown.

When the genetic model is known (e.g., recessive, dominant or additive model), the optimal Cochran-Armitage trend test with the appropriate choice of $x$ is more powerful than the chi-squared test with 2 df for testing an association. The genetic model is usually not known for complex diseases. In this situation, the maximum of three optimal tests (including the two extreme tests), $Z_{\text{MAX3}}$, is shown to be efficient robust compared to other available tests. In particular, $Z_{\text{MAX3}}$ is slightly more powerful than the chi-squared test with 2 df. Based on prior scientific knowledge, if the dominant model can be eliminated, then MERT, the maximum test, and the chi-squared test have roughly comparable power for a genetic model that ranges from recessive model to additive model and the allele frequency is not small. In this situation, the MERT and the chi-squared test are easier to apply than the maximum test and can be used by researchers. Otherwise, with current computational tools, $Z_{\text{MAX3}}$ is recommended.

### Acknowledgements

It is a pleasure to thank Prof. Javier Rojo for inviting us to participate in this important conference, in honor of Prof. Lehmann, and the members of the Department of Statistics at Rice University for their kind hospitality during it. We would also like to thank two referees for their useful comments and suggestions which improved our presentation.